\pdfoutput=1
\documentclass[11pt]{amsart}
\usepackage[a4paper,margin=27mm,headsep=6mm,footskip=11mm]{geometry}
\usepackage[T1]{fontenc}
\usepackage[utf8]{inputenc}
\usepackage{lmodern,microtype,mathtools,amssymb,booktabs,graphicx}
\usepackage{xcolor}
\usepackage{needspace,placeins}
\usepackage[colorlinks=true,linkcolor=blue!45!black,citecolor=blue!45!black,urlcolor=blue!45!black]{hyperref}
\usepackage[nameinlink,capitalize]{cleveref}
\crefname{figure}{Figure}{Figures}

\numberwithin{equation}{section}
\newtheorem{theorem}{Theorem}[section]
\newtheorem{proposition}[theorem]{Proposition}
\newtheorem{lemma}[theorem]{Lemma}
\newtheorem{corollary}[theorem]{Corollary}
\theoremstyle{remark}
\newtheorem{remark}[theorem]{Remark}
\AddToHook{env/proposition/begin}{\crefalias{theorem}{proposition}}
\AddToHook{env/lemma/begin}{\crefalias{theorem}{lemma}}
\AddToHook{env/corollary/begin}{\crefalias{theorem}{corollary}}
\AddToHook{env/remark/begin}{\crefalias{theorem}{remark}}
\newcommand{\C}{\mathbb C}
\newcommand{\R}{\mathbb R}
\newcommand{\Z}{\mathbb Z}

\newcommand{\conv}{\operatorname{conv}}

\newcommand{\M}{\mathcal M}
\newcommand{\A}{\mathcal A}
\newcommand{\doi}[1]{\href{https://doi.org/#1}{doi: #1}}
\title[Contact multiplicity and collinear dendrites]{Contact multiplicity and a sharp six-map bound\\ for collinear self-similar dendrites}
\author{Bernat Espigul\'e}
\address{Departament d'Inform\`atica, Matem\`atica Aplicada i Estad\'istica, Universitat de Girona, 17003 Girona, Catalonia, Spain}
\email{bernat@espigule.com}
\date{31 October 2026}
\hypersetup{pdftitle={Contact multiplicity and a sharp six-map bound for collinear self-similar dendrites},pdfauthor={Bernat Espigul\'e},pdfsubject={MSC2020: 28A80; 54F50, 30C15, 11C08},pdfkeywords={self-similar dendrite, collinear iterated function system, contact multiplicity, roots of integer polynomials, restricted power series, Jensen's formula}}

\begin{document}
\begin{abstract}
Let $K\subset\mathbb C$ be the attractor of $n\ge2$ similarities
$z\mapsto\lambda z+b_k$, $0<|\lambda|<1$, whose fixed points are evenly
spaced on a line, and suppose that $K$ does not lie in a line. We prove
that $K$ can be a dendrite only if $n\le6$, and that such dendrites
exist for every $2\le n\le6$; no separation or finite-intersection
hypothesis is imposed. Hence, for $n\ge7$, every connected $K$ of this
kind has interior points or separates the plane. Neighboring pieces of
$K$ meet in similar copies of a single contact set. In a dendrite exactly
one first-level piece meets this set, whereas for $n\ge7$ at least
$n-\lceil2\sqrt{n-1}\rceil\ge2$ do, a bound attained by explicit
quadratic parameters. The proof combines a planar translation
obstruction for replicated contact arcs with an estimate based on
Jensen's formula. After normalizing to $z\mapsto t+z/c$,
$t\in\{-n+1,-n+3,\ldots,n-1\}$, the connectedness locus is the closure,
outside the unit disk, of the roots of monic integer polynomials of
height at most $n-1$. The same estimate shows that, for $n\ge7$, every
such polynomial with a non-real root $c$ satisfying
$|c|^2\ge\lceil2\sqrt{n-1}\rceil+1$ has subleading coefficient of
absolute value at most $\lceil2\sqrt{n-1}\rceil$, and that this bound is
sharp. Accordingly, for $n\ge7$, at every non-real boundary point $c$ of
the locus with $|c|>1$, each expansion $c=\sum_{j\ge0}d_jc^{-j}$ with
integer digits $|d_j|\le n-1$ has $|d_0|\le\lceil2\sqrt{n-1}\rceil$.
\end{abstract}
\maketitle

\section{Introduction}\label{sec:intro}

Let $n\ge2$, let $\lambda\in\C$ with $0<|\lambda|<1$, and let $K\subset\C$
be the attractor~\cite{Hutchinson1981} of the similarities
\begin{equation}\label{eq:general-ifs}
 S_k(z)=\lambda z+b_0+k\delta\qquad(0\le k<n,\ \delta\ne0).
\end{equation}
The maps share their linear part, and their fixed points
$(b_0+k\delta)/(1-\lambda)$ are evenly spaced on a line; conversely,
every system of $n$ orientation-preserving similarities with a common
linear part and evenly spaced fixed points has this form. If $\lambda$
is real, the line through the fixed points is invariant and contains
$K$. If $\lambda$ is not real, $K$ lies in no line: a line $L\supseteq K$
would give $S_0(K)\subseteq L\cap S_0(L)$, a single point, whereas $S_0$
is injective and $K$ contains the distinct fixed points of $S_0$ and
$S_1$.

A \emph{dendrite} is a locally connected continuum containing no simple
closed curve; its points are joined by unique arcs
~\cite{CharatonikCharatonik1998,Nadler1992}. In the plane, a locally
connected continuum is a dendrite if and only if it has empty interior
and does not separate the plane. Indeed, a dendrite has empty interior,
since a disk contains a simple closed curve, and it is contractible, so
Alexander duality shows that it does not separate the plane; conversely,
if a
non-separating continuum contains a simple closed curve, then the
bounded complementary domain of that curve lies in the continuum. Planar
self-similar dendrites are thus exactly the self-similar continua that
are simply connected, in the sense of not separating the plane, and have
empty interior. For such a set, the identifications between symbolic
addresses must be compatible with the uniqueness of arcs. We ask how
many evenly spaced maps can produce one.

\begin{theorem}[Sharp six-map bound]\label{thm:main}
Let $K$ be the attractor of \eqref{eq:general-ifs}, and suppose that $K$
does not lie in a line. If $K$ is a dendrite, then $2\le n\le6$, and for
each $2\le n\le6$ there are such dendrites. In particular, for $n\ge7$
every connected $K$ of this form contains a simple closed curve, and
therefore has interior points or separates the plane.
\end{theorem}

Dendrites and other tree-like sets are basic models in analysis on
fractals. Kigami developed harmonic calculus on limits of networks with
dendrites as a principal application~\cite{Kigami1995}, and self-similar
dendrites such as Hata's tree-like set are standard examples of
postcritically finite self-similar sets~\cite{Kigami2001}. Dendrites
also arise as Julia sets: if every critical point of a polynomial is
strictly preperiodic, as for $z^2+i$, its Julia set is a dendrite
~\cite{DouadyHubbard1984,Milnor2006}. Julia sets of postcritically finite
rational maps are topologically self-similar~\cite{Kameyama2000}, and
Laplacians have been constructed on several quadratic Julia sets
~\cite{FlockStrichartz2012,RogersTeplyaev2010}. For binary systems with a
common multiplier and singleton overlap, the attractor carries a folded
inverse map; in the postcritically finite case, Ero\u{g}lu, Rohde, and
Solomyak prove a quasisymmetric conjugacy between this map and a
quadratic polynomial on its Julia set
~\cite[Theorem~1.1]{ErogluRohdeSolomyak2010}. \Cref{thm:main} shows that
evenly spaced systems with a common non-real multiplier produce such
tree-like attractors only when there are at most six maps.

\begin{figure}[!t]
\centering
\includegraphics[width=\textwidth,alt={Five fractal dendrites for n = 2 to 6, each a row of n colored, labeled pieces meeting at single circled points, above five magnified views of one contact.}]{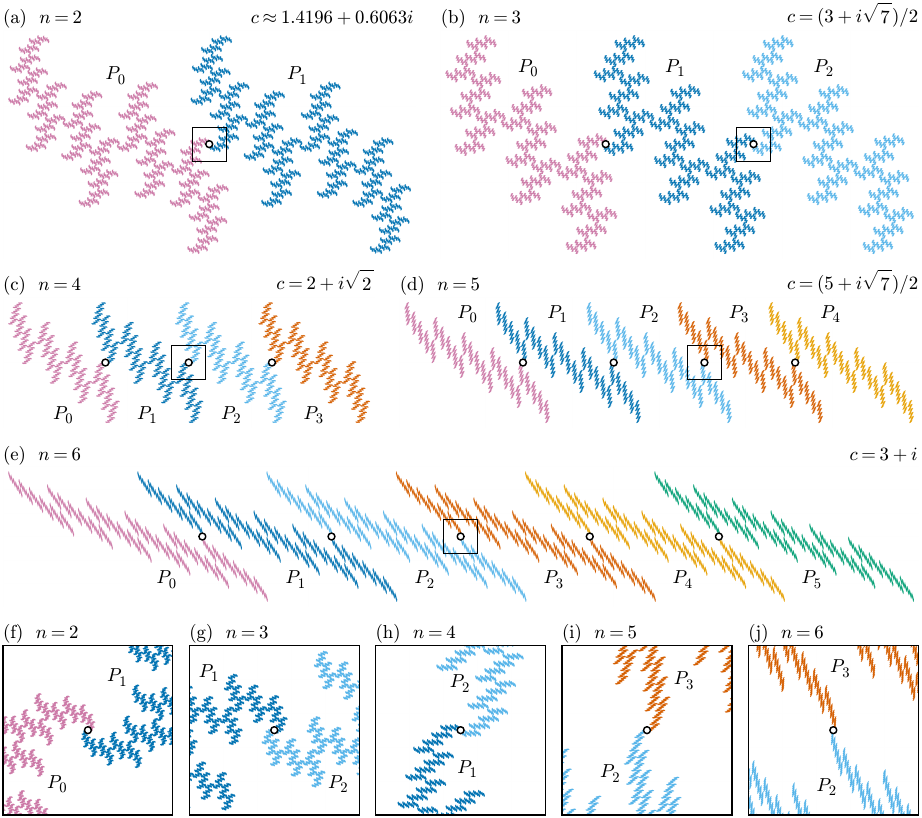}
\caption{Non-real dendrites at all five admissible alphabet sizes, with
the exact parameters in \eqref{eq:examples}. Colors and labels identify
the first-level pieces $P_0,\dots,P_{n-1}$ of \eqref{eq:pieces}, in
increasing digit order; white circles mark the contact points
$P_i\cap P_{i+1}=\{a_i+1\}$. Panels (a)--(e) are drawn at different
scales, and panels (f)--(j) magnify the squares in (a)--(e) five times.
The point samples illustrate the sets; exact intersection certificates and proofs are given in \cref{sec:examples}.}
\label{fig:examples}
\end{figure}

The affine map $z\mapsto2(z-w)/\delta$, with
$w=(b_0+(n-1)\delta/2)/(1-\lambda)$, conjugates \eqref{eq:general-ifs} to
the system $f_t(z)=t+c^{-1}z$, $t\in A_n$, where $c=1/\lambda$ and
$A_n=\{-n+1,-n+3,\ldots,n-1\}$. It carries $K$ onto
\begin{equation}\label{eq:attractor}
 E(c,n)=\left\{\sum_{j\ge0}t_jc^{-j}:t_j\in A_n\right\},
 \qquad |c|>1,
\end{equation}
and $K$ lies in a line exactly when $c$ is real. We work in this
normalization throughout. Write
$\M_n=\{c:|c|>1,\ E(c,n)\text{ is connected}\}$ for its
connectedness locus. In these terms, \cref{thm:main} states that $E(c,n)$
is a dendrite for some non-real $c$ exactly when $2\le n\le6$;
\cref{fig:examples} shows one for each such $n$.

Barnsley and Harrington initiated the parameter-space study of the binary
family~\cite{BarnsleyHarrington1985}. They considered the disconnectedness
set for a pair of similarities with a common complex multiplier. The
connectedness convention and reciprocal change $\lambda=1/c$ give
$\M_2$. Bousch proved connectedness, local connectedness, and H\"older
path connectedness of the binary locus~\cite{Bousch1993}.
His thesis also treats the consecutive alphabets $\{0,1,\ldots,n-1\}$
for every $n\ge2$, proving that connectedness is equivalent to overlap
between neighboring first-level pieces
~\cite[Chapter~1, p.~41]{Bousch1992}. An affine change of coordinates
identifies these systems with \eqref{eq:attractor}.
Bandt~\cite{Bandt2002} and Solomyak and Xu~\cite{SolomyakXu2003}
developed the study of the binary locus, Shmerkin and Solomyak related
zeros of $\{-1,0,1\}$ power series to connectedness loci of self-affine
sets~\cite{ShmerkinSolomyak2006}, and Calegari, Koch, and Walker proved
Bandt's conjecture on the binary locus~\cite{CalegariKochWalker2017}.
Espigul\'e, Juher, and Salda\~na~\cite{EspiguleJuherSaldana2024}
investigated the collinear loci uniformly in $n$. When $c$ is a non-real
quadratic integer with $|c|^2=n$, the set $E(c,n)$ is a planar self-affine tile
with a collinear digit set; the connectedness, disklikeness, and
topology of such tiles were studied by Kirat and Lau~\cite{KiratLau2000},
Leung and Lau~\cite{LeungLau2007}, and Akiyama, Loridant, and
Thuswaldner~\cite{AkiyamaLoridantThuswaldner2021}. Bandt and Hung studied
the analogous family whose fixed points are the vertices of a regular
polygon~\cite{BandtHung2008}.

\Needspace{8\baselineskip}
The non-real restriction is essential: $E(2,n)=[-2(n-1),2(n-1)]$
is an interval, hence a dendrite, for every $n\ge2$.
Connected self-similar sets are locally connected~\cite{Hata1985}, so
the exclusion amounts to forcing a simple closed curve. Finite overlaps
are relevant to this problem: Bandt and Rao prove that finite first-level
intersections in a connected planar self-similar set imply the open set
condition~\cite{BandtRao2007}. Self-similar dendrites can nevertheless have
nontrivial subdendrites as intersections, even under that condition
~\cite{AllabergenovaSamuelTetenov2024}. Our argument uses the equal spacing
of the digits and applies without a separation or finite-intersection
hypothesis.

The exclusion follows from a quantitative statement about contacts. Write
\begin{equation}\label{eq:pieces}
 \begin{gathered}
 K=E(c,n),\qquad a_i=-n+1+2i,\qquad P_i=a_i+c^{-1}K\quad(0\le i<n),\\
 Q=K\cap(2c+K),
 \end{gathered}
\end{equation}
so that $P_i\cap P_{i+1}=a_i+c^{-1}Q$ for $0\le i<n-1$.
Thus $Q$ is the common normalized contact between neighboring pieces.
Define
\begin{equation}\label{eq:k}
 k(c,n)=\#\{i:P_i\cap(2c+K)\ne\varnothing\}.
\end{equation}
Since $P_i\subseteq K$, equivalently
$k(c,n)=\#\{i:P_i\cap Q\ne\varnothing\}$. This counts participating
pieces, not points or connected components of $Q$.

\Needspace{11\baselineskip}
\begin{theorem}[Contact multiplicity]\label{thm:multiplicity}
Let $c\in\M_n\setminus\R$.
\begin{enumerate}
\item If $E(c,n)$ is a dendrite, then exactly one of
$P_{n-1}\cap(2c+P_0)$ and $P_0\cap(2c+P_{n-1})$ is nonempty, all other
sets $P_i\cap(2c+P_j)$ are empty, and $k(c,n)=1$.
\item If $n\ge7$, then
\begin{equation}\label{eq:uniform-count}
 k(c,n)\ge k_n:=n-\lceil2\sqrt{n-1}\rceil\ge2.
\end{equation}
\end{enumerate}
For every $n\ge7$, equality is attained at a non-real quadratic
parameter. Thus $k_n$ is the exact minimum of $k(\cdot,n)$ over
$\M_n\setminus\R$.
The sequence $(k_n)_{n\ge7}$ is nondecreasing and
${k_n=n-2\sqrt n+O(1)}$.
\end{theorem}

For example, $k_7=k_8=2$, $k_9=3$, $k_{10}=4$,
$k_{20}=11$, and $k_{100}=80$. On the imaginary axis every piece participates:
$k(ir,n)=n$ whenever $1<r\le\sqrt n$ (\cref{prop:imaginary}).

The proof has two independent parts. Equal digit offsets replicate
contacts across consecutive pieces. If two such contacts occur inside a
dendrite, their unique arcs produce a forbidden translation
parallelogram. Analytically, a sufficiently large first digit forces a
small real zero of the corresponding power series. Jensen's formula then
excludes a non-real zero at the required modulus, in the spirit of the
estimates of Beaucoup, Borwein, Boyd, and Pinner~\cite{BBBP1998}. Connectivity at a
smaller alphabet handles the complementary range of moduli.

The same estimate constrains roots of integer polynomials. For $n=2$,
the description of the connectedness locus through zeros of power series
and polynomials with coefficients in $\{-1,0,1\}$ is classical
~\cite{CalegariKochWalker2017,ShmerkinSolomyak2006}. For every $n$,
$\M_n$ is the closure, outside the closed unit disk, of the set of roots
of monic integer polynomials of height at most $n-1$, that is, with all
coefficients in $\{-n+1,\ldots,n-1\}$
~\cite[Proposition~2.3]{EspiguleJuher2026Layers}. The closure in $\C$ of
the roots of all integer polynomials of height at most $n-1$, monic or
not, is $\M_n\cup\{|z|=1\}\cup\{0\}\cup\{1/c:c\in\M_n\}$
~\cite[Proposition~1.1]{EspiguleJuher2026Layers}; in particular,
non-monic polynomials add nothing outside the closed unit disk. A root $c$ of
$z^k+a_1z^{k-1}+\cdots+a_k$ satisfies
$c=-a_1-a_2c^{-1}-\cdots-a_kc^{1-k}$, a finite digit expansion of the kind
in \cref{lem:connectivity}, whose first digit is $-a_1$.

\begin{corollary}[Subleading coefficients of non-real roots]\label{thm:roots}
Let $n\ge7$, and let $c\notin\R$ with
$|c|^2\ge\lceil2\sqrt{n-1}\rceil+1$ be a root of a monic polynomial
$z^k+a_1z^{k-1}+\cdots+a_k$ with integer coefficients of absolute value
at most $n-1$. Then $|a_1|\le\lceil2\sqrt{n-1}\rceil$. The bound is
sharp: for every $n\ge7$, some such root has $|c|^2>n$ and
$|a_1|=\lceil2\sqrt{n-1}\rceil$.
\end{corollary}

The proof is given in \cref{sec:piece-loci}, after
\cref{thm:exterior-labels}, because the sharpness examples use the
quadratic parameters introduced there.

Since $\lceil2\sqrt{n-1}\rceil+1\le n-1$ for $n\ge7$, the corollary applies
in particular to every non-real root with $|c|\ge\sqrt n$; for real
coefficients of absolute value at most $n-1$, the same argument gives
$|a_1|<\lceil2\sqrt{n-1}\rceil+1$. For such roots the extreme values
$a_1=\pm(n-1)$ are excluded from $n=7$ on, whereas for $3\le n\le6$ the
construction in the proof of \cref{thm:roots} gives non-real roots with
$|c|^2>n$ and $|a_1|=n-1$. Since
$\{1<|c|<\sqrt n\}\subseteq\M_n$ by \cref{lem:connectivity} and
$|c|<1+\sqrt{n-1}$ on $\M_n\setminus\R$ by \cref{prop:piece-strips}, the
non-real points of $\partial\M_n$ with $|c|>1$, which belong to $\M_n$
because $\M_n$ is closed in $\{|c|>1\}$
~\cite[Proposition~2.3]{EspiguleJuher2026Layers}, lie in the annulus
$\sqrt n\le|c|<1+\sqrt{n-1}$. There, by \cref{thm:exterior-labels},
every expansion $c=\sum_{j\ge0}d_jc^{-j}$ with integer digits
$|d_j|\le n-1$ has $|d_0|\le\lceil2\sqrt{n-1}\rceil<n-1$ once $n\ge7$.
Geometrically, the non-real parts of the parameter regions
selected by large first digits, and even of their convex relaxations,
recede inside the connectedness disk from $n=7$ on;
\cref{fig:piece-cells,fig:piece-atlas} display these vanishing fronts. In
\cref{sec:piece-loci} we determine exactly which first digits occur at
non-real parameters outside the connectedness disk, give support-function
bounds for the corresponding parameter sets, and use a quadratic family
to attain the minimum in \cref{thm:multiplicity}. \Cref{fig:contact}
shows such a parameter, $c=3+i\sqrt2$ with $n=7$ and $k(c,7)=2$: the
pieces $P_5$ and $P_6$ meet $2c+P_0$ and $2c+P_1$, respectively, while
$P_6$ does not meet $2c+P_0$. In a dendrite, by contrast, the only
contact is between $P_{n-1}$ and $2c+P_0$, or $P_0$ and $2c+P_{n-1}$
(\cref{thm:multiplicity}(1)).

The examples in \cref{fig:examples} establish sharpness at the following
exact parameters:
\begin{equation}\label{eq:examples}
 \begin{cases}
 n=2\colon & c\text{ is the root in the upper half-plane of }c^3-2c^2+2=0,\\
 3\le n\le6\colon & c=\dfrac{n+i\sqrt{8(n-1)-n^2}}2.
 \end{cases}
\end{equation}
Up to affine normalization and complex conjugation, the binary example
is the one in~\cite[Figure~3]{ErogluRohdeSolomyak2010}.
In each case neighboring pieces meet at one point, non-neighboring pieces
are disjoint, and the first-level intersection graph is a path.
\Cref{sec:examples} proves these assertions by finite exact
inverse-state calculations. The computations are included in full and
can also be reproduced by the accompanying verification program.

\subsection*{Relation to other work}
The present classification concerns the possible alphabet sizes; a
classification of all dendrite parameters, even at a fixed size, is a
finer problem.

Recent work on the same collinear family concerns the geometry of
$\M_n$. Espigul\'e and Juher construct semialgebraic interior layers
using finite capture and prove uniform bounds on enclosure-tree branching
~\cite{EspiguleJuher2026Layers}. A separate paper proves that every
non-real point of $\M_n$ is a limit of interior points
~\cite{Espigule2026RegularClosedness}, and infinitely many holes of the
loci are exhibited in~\cite{Espigule2026Holes}. Here the object of study is the
topology of the individual attractor. Our contact count is linked to the
inverse search through the replication formula in
\cref{lem:replication}; \cref{rem:branching-contacts} explains the
distinction between the number of admissible offsets and the number of
pieces they connect. The normalization follows the thesis
~\cite{Espigule2026Thesis}, whose Chapter~3 and Appendix~C announce the
dendrite classification. The present paper gives the proofs and the
sharp minimum contact count.

\section{Connectivity and replicated contacts}\label{sec:contacts}

We first recall Bousch's connectedness criterion for consecutive collinear
digits~\cite[Chapter~1, p.~41]{Bousch1992}. We include a proof together
with the separation consequence and the classical connectedness-disk
estimate~\cite[Lemma~2.5]{Solomyak1998}, in the form of
~\cite[Proposition~5(i)]{EspiguleJuherSaldana2024}, needed below.
For $p\ge2$, put $D_p=\{-p+1,-p+2,\ldots,p-1\}$. Digit subtraction gives
\begin{equation}\label{eq:balanced}
 \frac12\bigl(E(c,p)-E(c,p)\bigr)
 =\left\{\sum_{j\ge0}d_jc^{-j}:d_j\in D_p\right\}.
\end{equation}
Every digit in $D_p$ can be lifted to a pair of digits from $A_p$;
thus this equality holds in both directions.

\Needspace{9\baselineskip}
\begin{lemma}[Neighboring contact and separation]\label{lem:connectivity}
For $|c|>1$ and $p\ge2$, the following are equivalent:
\begin{enumerate}
\item $E(c,p)$ is connected;
\item two distinct first-level pieces intersect;
\item every neighboring pair of first-level pieces intersects;
\item $c=\sum_{j\ge0}d_jc^{-j}$ for some $d_j\in D_p$.
\end{enumerate}
If $E(c,p)$ is disconnected, its first-level pieces are pairwise
disjoint. Moreover,
\begin{equation}\label{eq:inner-disk}
 1<|c|<\sqrt p\quad\Longrightarrow\quad c\in\M_p.
\end{equation}
\end{lemma}

\begin{proof}
We use the elementary planar crossing fact that a left-to-right continuum
and a bottom-to-top continuum in a rectangle must meet; see, for example,
~\cite[Lemma~7]{PascolettiZanolin2013}. The argument below is Bousch's
\emph{lemme de coin\c{c}age}~\cite[Chapter~1, p.~28]{Bousch1992}.
The same overlap criterion can be obtained from Hopf's horizontal-chord
theorem~\cite{Hopf1936}; see~\cite[Proposition~2.2 and Lemma~A.1]{EspiguleJuher2026Layers}.
Let $L$ be a continuum and suppose two translates $a_i+L$ and $a_j+L$
intersect. Choose the smallest positive index gap $s$ giving such an
intersection. Suppose that $L$ has positive vertical extent and $s\ge2$,
and let $R$ be the bounding rectangle of the continuum
$(a_i+L)\cup(a_{i+s}+L)$. This continuum meets both vertical sides of
$R$. The intermediate translate $a_{i+1}+L$ lies in $R$ and meets both
horizontal sides, because all translates have the vertical extent of
$L$. An intersection gives a gap $1$ or $s-1$, a contradiction.
Thus $s=1$. If $L$ has zero vertical extent, it is a horizontal interval
and the same conclusion follows from its length. Translation invariance
then makes all neighboring translates intersect.

To apply this without assuming the attractor connected, let
$K_p=E(c,p)$, $C_0=\conv K_p$, and
\[
 C_{r+1}=\bigcup_{t\in A_p}(t+c^{-1}C_r).
\]
These compact sets are nested, contain $K_p$, and converge to $K_p$.
Suppose two actual first-level pieces intersect. Starting with $C_0$,
inductively assume $C_r$ connected. The translates of $c^{-1}C_r$
inherit that intersection, so the preceding crossing argument shows that
all neighboring translates meet. Hence $C_{r+1}$ is connected.
Their nested intersection $K_p$ is therefore connected. Conversely, a
connected finite union of nonempty compact pieces cannot have those
pieces pairwise disjoint. Applying the crossing argument directly to the
connected set $c^{-1}K_p$ proves the equivalence of the first three
conditions. The fourth follows from \eqref{eq:balanced}, since neighboring
contact is equivalent to $2c\in K_p-K_p$.

If $K_p$ is disconnected, the established equivalence gives positive
distance between every pair of first-level pieces. For sufficiently small
$\varepsilon>0$, let $O$ be the open $\varepsilon$-neighborhood of $K_p$.
The images $f_t(O)$ are pairwise disjoint, and each lies in $O$, since
its points have distance less than $\varepsilon/|c|$ from $K_p$.
Taking planar area, which is finite and positive on $O$, gives
$p|c|^{-2}\operatorname{area}(O)\le\operatorname{area}(O)$.
Thus $|c|^2\ge p$, proving \eqref{eq:inner-disk}.
\end{proof}

For a connected $K$, each consecutive block $P_i\cup\cdots\cup P_j$
is consequently a continuum. This fact will control the arcs used in the
topological argument.

Set $h=n-1$, the height bound in \cref{thm:roots}, so that
$D_n=\{-h,\ldots,h\}$, and introduce the admissible offsets
\begin{equation}\label{eq:offsets}
 \A(c,n)=\left\{m\in\Z:|m|\le h,\quad
 c=m+\sum_{j\ge1}d_jc^{-j}\text{ for some }d_j\in D_n\right\}.
\end{equation}
By \cref{lem:connectivity}, this set is nonempty whenever $c\in\M_n$.

The replication of overlaps under equal digit differences already appears
in Bousch's translation symmetries of the overlap graph
~\cite[Chapter~1, p.~38]{Bousch1992}. The next lemma records this mechanism
for the contact $K\cap(2c+K)$ and counts the resulting pairs.

\begin{lemma}[Contact replication]\label{lem:replication}
For $0\le i,j<n$,
\begin{equation}\label{eq:pair-offset}
 P_i\cap(2c+P_j)\ne\varnothing
 \quad\Longleftrightarrow\quad i-j\in\A(c,n).
\end{equation}
Each $m\in\A(c,n)$ gives exactly $n-|m|$ intersecting ordered pairs,
whose intersection sets are successive translates by $2$. In particular,
\begin{equation}\label{eq:count-formula}
 J(c,n):=\#\{(i,j):P_i\cap(2c+P_j)\ne\varnothing\}
 =\sum_{m\in\A(c,n)}(n-|m|),
\end{equation}
and
\begin{equation}\label{eq:k-lower}
 k(c,n)\ge n-\min_{m\in\A(c,n)}|m|.
\end{equation}
\end{lemma}

\begin{proof}
An equality $a_i+u/c=2c+a_j+v/c$, with $u,v\in K$, is equivalent to
$u-v=2c(c-(i-j))$. Using \eqref{eq:balanced} gives
\eqref{eq:pair-offset}. For a fixed $m$, there are exactly $n-|m|$
pairs with $i-j=m$. Incrementing both indices by one translates the
intersection by $2$. Their left indices are distinct, which proves
\eqref{eq:k-lower}; summation gives \eqref{eq:count-formula}.
\end{proof}

\section{The topological obstruction}\label{sec:topology}

We use standard properties of dendrites: subcontinua are dendrites;
nonempty intersections of subcontinua are connected; each pair of points
is joined by a unique arc; and a dendrite is contractible
~\cite{CharatonikCharatonik1998}.

\begin{lemma}[Translation parallelogram]\label{lem:parallelogram}
Let $T$ be a planar dendrite and let $u,v\in\C$ be linearly independent
over $\R$. There are no paths $\alpha$ from $p$ to $p+u$ and $\beta$
from $p$ to $p+v$ such that
\[
 \alpha,\quad\alpha+v,\quad\beta,\quad\beta+u
\]
all lie in $T$.
\end{lemma}

\begin{proof}
These four paths define a continuous map on the boundary of the unit
square, using identical parametrizations on opposite sides. Contractibility
of $T$ extends it to a continuous map $F\colon[0,1]^2\to T$. Extend $F$ over
the plane by
\[
 H(r+x,s+y)=F(x,y)+ru+sv
 \quad(r,s\in\Z,\ 0\le x,y\le1).
\]
The prescribed boundary maps agree on shared edges, so $H$ is well-defined
and continuous. For the invertible real-linear map $L(x,y)=xu+yv$,
the function $H-L$ is periodic and bounded. For any $z\in\C$, the map
\[
 x\longmapsto L^{-1}z-L^{-1}\bigl(H(x)-Lx\bigr)
\]
takes a sufficiently large closed disk into itself. Brouwer's fixed-point
theorem shows that $H(x)=z$ has a solution. Thus $H$ is onto.

But its image is contained in
$\bigcup_{r,s\in\Z}(T+ru+sv)$. A planar dendrite has empty interior,
since a disk contains a simple closed curve; its translates are compact
and nowhere dense. Baire's theorem rules out this covering of the plane.
\end{proof}

\begin{proposition}[One participating piece]\label{prop:one-piece}
If $c\notin\R$ and $K=E(c,n)$ is a dendrite, then precisely one of
\[
 P_{n-1}\cap(2c+P_0),\qquad
 P_0\cap(2c+P_{n-1})
\]
is nonempty, and all other ordered cross-contacts are empty. Consequently
$k(c,n)=1$.
\end{proposition}

\begin{proof}
The contact $Q$ is nonempty by \cref{lem:connectivity}. For adjacent digits,
the union of the corresponding pieces is a subcontinuum of $K$, and is
similar to $K\cup(2c+K)$. The latter union is therefore a dendrite. Its
subcontinua $K$ and $2c+K$ have connected intersection $Q$.

Suppose $m\in\A(c,n)$ and $|m|\le n-2$. Choose $i,j$ with $i-j=m$
and $i+1,j+1<n$, and take
$x\in P_i\cap(2c+P_j)$. Replication gives
$x+2\in P_{i+1}\cap(2c+P_{j+1})$. The unique arc $\gamma$ from $x$
to $x+2$ lies in $Q$, so both $\gamma$ and $\gamma-2c$ lie in $K$.
The unique arc $\delta$ from $x-2c$ to $x$ lies in the consecutive block
\[
 B=P_{\min(i,j)}\cup\cdots\cup P_{\max(i,j)},
\]
because its endpoints lie in $P_j$ and $P_i$ and $B$ is a subcontinuum.
The shifted block $B+2$ also lies in $K$, giving $\delta+2\subseteq K$.
Now apply \cref{lem:parallelogram} with
\[
 p=x-2c,\quad u=2,\quad v=2c,\quad
 \alpha=\gamma-2c,\quad\beta=\delta.
\]
The required vectors are independent because $c$ is non-real. This
contradiction excludes every nonextreme offset.

Let $Q_+=P_{n-1}\cap(2c+P_0)$ and $Q_-=P_0\cap(2c+P_{n-1})$.
We have $Q=Q_+\cup Q_-$. If both sets were nonempty, connectedness of $Q$
would force $Q_+\cap Q_-\ne\varnothing$. A point of their intersection
would lie in $P_0\cap(2c+P_0)$, producing the already excluded offset
$0$. Thus exactly one extreme contact occurs. A point of $Q$ in any
other left piece would produce an additional cross-contact, so exactly
one left piece participates.
\end{proof}

\section{The sharp uniform contact bound}\label{sec:jensen}

Zeros of power series $1+\sum_{i\ge1}a_iz^i$ with real coefficients in
$[-g,g]$ were studied by Beaucoup, Borwein, Boyd, and Pinner
~\cite{BBBP1998,BBBP1998b}. Among other results, they show that for
$g\ge3+2\sqrt2$ the smallest disk containing two zeros has radius
$(\sqrt g+1)^{-1}$, the smallest possible double real zero, and they
compute the smallest circle containing a real zero and a pair of
conjugate zeros of a given argument~\cite{BBBP1998}. The next lemma
adapts these ideas to a prescribed large first coefficient, using only a
bound on the remaining coefficients, which need not be integers. Its
hypothesis $h\ge6$ is the integer form of $1+2\sqrt h\le h$, that is,
$h\ge3+2\sqrt2$. For $m=1+2\sqrt h$, the series
$1-mz+h\sum_{j\ge2}z^j$ used for sharpness below has the double zero
$(1+\sqrt h)^{-1}$.

\begin{lemma}[First-coefficient exclusion]\label{lem:jensen}
Let $h\ge6$ be an integer and $q=\lceil2\sqrt h\rceil+1$. If
\[
 F(z)=1-mz-\sum_{j\ge2}b_jz^j,
 \qquad q\le |m|\le h,\quad m,b_j\in\R,\quad |b_j|\le h,
\]
then every zero of $F$ in $|z|\le q^{-1/2}$ is real.
\end{lemma}

\begin{proof}
Replacing $F(z)$ by $F(-z)$ if necessary, take $m\ge q$.
Notice that $6\le q\le h$, since $4<2\sqrt h\le h-1$ for $h\ge6$, and
$m\ge1+2\sqrt h$. For $0<t<1$,
\[
 F(t)\le U(t):=1-mt+\frac{ht^2}{1-t}
 =\frac{1-(m+1)t+(m+h)t^2}{1-t}.
\]
The numerator has smaller positive root
\begin{equation}\label{eq:majorant-root}
 t_0=\frac{2}{m+1+\sqrt{(m-1)^2-4h}}
 \le\frac1{1+\sqrt h},\qquad t_0<\frac2m.
\end{equation}
Since $F(0)=1$ and $F(t_0)\le0$, there is a positive real zero
$\alpha\le t_0$.

Suppose it also has a non-real zero $\lambda$ with $|\lambda|^2\le1/q$.
Its conjugate is another zero. Set $\rho=1/\sqrt2$; all three distinct
zeros lie in $|z|<\rho$. Jensen's formula, followed by the
arithmetic--geometric mean inequality and Parseval's identity, gives
\begin{align}
 \log\frac{\rho^3}{\alpha|\lambda|^2}
 &\le\frac1{2\pi}\int_0^{2\pi}\log|F(\rho e^{it})|\,dt\notag\\
 &\le\frac12\log\left(\frac1{2\pi}
       \int_0^{2\pi}|F(\rho e^{it})|^2\,dt\right)
 \le\frac12\log\left(1+\frac{m^2}{2}+\frac{h^2}{2}\right).
 \label{eq:jensen-estimate}
\end{align}
If a zero lies on $|z|=\rho$, obtain the same estimate by a limit through
radii avoiding the zeros. On the other hand,
\[
 \frac{\rho^3}{\alpha|\lambda|^2}>
 \frac{mq}{4\sqrt2}>
 \sqrt{1+\frac{m^2}{2}+\frac{h^2}{2}}.
\]
For the last inequality, it suffices that
$m^2(q^2-16)>16h^2+32$. Its left side is at least $q^2(q^2-16)$.
For $h=6$, subtracting the right side gives $112$. For $h\ge7$,
using $q\ge2\sqrt h+1$ and $\sqrt h>5/2$ gives
\[
 q^2(q^2-16)-16h^2-32
 \ge(32h-56)\sqrt h-40h-47
 >40h-187>0.
\]
This contradicts \eqref{eq:jensen-estimate} and proves the lemma.
\end{proof}

The middle expression in the last display is positive for every
$h\ge6$; it equals $46.1\ldots$ at $h=6$ and increases with $h$. Hence the
same proof applies with $q$ replaced by the real number $1+2\sqrt h$,
both in the bound on $|m|$ and in the radius. This threshold cannot be
lowered: for real $m$ slightly below $1+2\sqrt h$, the series
$1-mz+h\sum_{j\ge2}z^j$ has a pair of non-real zeros of modulus
$(m+h)^{-1/2}<(1+2\sqrt h)^{-1/2}$, which merge at $m=1+2\sqrt h$ into
the double zero $(1+\sqrt h)^{-1}$.

\begin{proof}[Proof of \cref{thm:multiplicity}]
Part~(1) is \cref{prop:one-piece}. For part~(2), write
$h=n-1\ge6$ and $q=\lceil2\sqrt h\rceil+1\le h$.
If $|c|^2<q$, then $E(c,q)$ is connected by
\cref{lem:connectivity}. It supplies an expansion of $c$ with all digits
in $D_q$, hence an admissible first offset $m$ with $|m|\le q-1$ for
the $n$-map system. By \eqref{eq:k-lower}, $k(c,n)\ge n-q+1$.

If $|c|^2\ge q$, take any $m\in\A(c,n)$. Its expansion in
\eqref{eq:offsets}, multiplied by $c^{-1}$, makes $\lambda=c^{-1}$ a
zero of a series of the form in \cref{lem:jensen}, with tail coefficients
bounded by $h$. As $c$ is non-real, that lemma excludes $|m|\ge q$.
Thus every admissible offset has $|m|\le q-1$, and
\eqref{eq:k-lower} again gives $k(c,n)\ge n-q+1$.

The inequality $q\le h$ gives $k_n\ge2$. The increase
$2\sqrt{h+1}-2\sqrt h$ is less than $1$, so the corresponding ceiling
increases by at most $1$. Hence $k_n$ is nondecreasing. Its asymptotic
form follows directly from the ceiling formula. Attainment is proved in
\cref{prop:attainment} below.
\end{proof}

\begin{proof}[Proof of the exclusion in \cref{thm:main}]
By the conjugation in \cref{sec:intro}, we may assume that $K=E(c,n)$
with $c$ non-real. For $n\ge7$ the two conclusions of
\cref{thm:multiplicity} are incompatible with a dendrite. A connected
self-similar set is locally connected~\cite{Hata1985}, so a connected
attractor that is not a dendrite contains a simple closed curve $\Gamma$.
If $K$ contains the bounded complementary domain of $\Gamma$, it has
interior points. Otherwise some point of that domain lies outside $K$,
and its component in $\C\setminus K$ lies inside $\Gamma$; hence $K$
separates the plane.
\end{proof}

\begin{remark}[The role of Jensen's formula]\label{rem:elementary}
For $n\ge10$ the exclusion also follows from elementary estimates. Let
$c$ be a non-real dendrite parameter and $h=n-1$. By
\cref{prop:one-piece}, $\A(c,n)$ is $\{h\}$ or $\{-h\}$; replacing $c$ by
$-c$ or $\bar c$, which leaves the attractor unchanged up to a
reflection, we may assume that $\A(c,n)=\{h\}$ and $\operatorname{Im}c>0$.
Then \eqref{eq:radial} with $k(c,n)=1$ gives $r^2\ge h$ for $r=|c|$, and
\cref{prop:piece-strips} gives $r<1+\sqrt h$ and
\[
 h\le2\operatorname{Re}c+\frac h{(r-1)^2}\le2r+\frac h{(r-1)^2}.
\]
The convex function
$f(r)=2r+h/(r-1)^2-h$ satisfies
$f(1+\sqrt h)=-(\sqrt h-3)(\sqrt h+1)\le0$ when $h\ge9$. Also
$f(\sqrt h)<0$ for $h\ge9$: with $s=\sqrt h$, the value
$f(\sqrt h)=2s+s^2/(s-1)^2-s^2$ decreases for $s>1$ and equals $-3/4$
at $s=3$. Hence $f<0$ on $[\sqrt h,1+\sqrt h)$, a contradiction.
\Cref{lem:jensen} is therefore
needed only for $7\le n\le9$ and for the sharp count in
\cref{thm:multiplicity}.
\end{remark}

\section{Dependence on the parameter}\label{sec:parameter}

The uniform estimate can be supplemented by elementary bounds that retain
the modulus of $c$ or the alphabet size.

\begin{proposition}[Radial bound]\label{prop:parameter}
Let $c\in\M_n\setminus\R$ and $r=|c|$. Then
\begin{equation}\label{eq:radial}
 k(c,n)\ge\max\{1,n-\lfloor r^2\rfloor\}.
\end{equation}
For $n\ge7$ one may take the maximum of this bound and $k_n$.
\end{proposition}

\begin{proof}
If $\lfloor r^2\rfloor+1\le n$, apply
\cref{lem:connectivity} to $p=\lfloor r^2\rfloor+1>r^2$ and then
\cref{lem:replication}. Otherwise the bound is just $k\ge1$.
\end{proof}

\begin{proposition}[Increasing the alphabet]\label{prop:alphabet}
Fix $c\in\M_n$. Then
\begin{equation}\label{eq:alphabet-growth}
 k(c,n+1)\ge k(c,n)+1.
\end{equation}
In particular, writing $\nu(c)=\min\{p\ge2:c\in\M_p\}$, one has
$k(c,n)\ge n-\nu(c)+1$ for $n\ge\nu(c)$. A non-real dendrite
$E(c,n)$ can occur only at the first connected alphabet size: $n=\nu(c)$.
\end{proposition}

\begin{proof}
The nested balanced digit alphabets give
$\A(c,n)\subseteq\A(c,n+1)$. An offset $m$ contributes exactly the
left indices
\[
 I_m^{(n)}=
 \begin{cases}
 \{m,\ldots,n-1\},&m\ge0,\\
 \{0,\ldots,n-1+m\},&m<0.
 \end{cases}
\]
Their union over $\A(c,n)$ is precisely the participating-index set.
Every old interval is contained in its counterpart at alphabet size $n+1$.
If there is a nonnegative offset, the new index $n$ is added. If all
offsets are negative, the union is an initial interval, and its endpoint
increases by at least one. This proves \eqref{eq:alphabet-growth}.
Iteration gives the bound involving $\nu(c)$; the dendrite assertion
then follows from \cref{prop:one-piece}.
\end{proof}

\begin{proposition}[Imaginary-axis count]\label{prop:imaginary}
For $r>1$, $E(ir,n)$ is connected if and only if $r^2\le n$.
Throughout that range it is a nondegenerate rectangle and
$k(ir,n)=n$.
\end{proposition}

\begin{proof}
Splitting the digit expansion into even and odd positions gives
\[
 E(ir,n)=S-iS/r,\qquad
 S=\left\{\sum_{j\ge0}t_j(-r^{-2})^j:t_j\in A_n\right\}.
\]
The convex hull of $S$ is $[-L,L]$, with
$L=hr^2/(r^2-1)$. Its $n$ first-level intervals cover this hull
precisely when $L/r^2\ge1$, or $r^2\le n$. In this case their invariant
set is the whole interval. If $r^2>n$, those intervals are pairwise
disjoint, so $S$ and hence $E(ir,n)$ are disconnected.
In the connected case, putting $\eta=hr/(r^2-1)$ yields
\[
 K=[-r\eta,r\eta]+i[-\eta,\eta],\qquad
 P_i=[a_i-\eta/r,a_i+\eta/r]+i[-\eta,\eta].
\]
Since $\eta\ge r$, every $P_i$ meets $2ir+K$. This includes equality
$r^2=n$, where $K\cap(2ir+K)=[-r\eta,r\eta]+ir$ is a horizontal segment.
\end{proof}

\begin{remark}[Branching and contact multiplicity]\label{rem:branching-contacts}
Put $N=2n-1$ and $g_t(z)=c(z-t)$. For $m\in D_n$, \eqref{eq:offsets} gives
\[
 m\in\A(c,n)\quad\Longleftrightarrow\quad
 g_{2m}(2c)\in E(c,N).
\]
Thus admissible offsets are precisely the first inverse branches of the
marked point $2c$ that remain in the difference attractor. When
$\A(c,n)\ne\varnothing$, also $2c\in E(c,N)$ by \cref{lem:connectivity}.
Since $E(c,N)$ lies in the canonical enclosure of
~\cite[Section~6.1]{EspiguleJuher2026Layers}, each digit $2m$ with
$m\in\A(c,n)$ then labels a child of the root $2c$ in the exact
enclosure tree of~\cite[Definition~6.3]{EspiguleJuher2026Layers}.
Theorem~E of that paper therefore gives $\#\A(c,n)\le4$ for non-real $c$
with $|c|\ge\sqrt n$, and $\#\A(c,n)\le3$ if $n=2$ or $n\ge25$. For
$n\ge10$ and $n\ge26$, respectively, these two bounds also follow from
\cref{prop:piece-strips}: when $|c|\ge\sqrt n$, the admissible offsets lie
in a closed interval of length at most $2+4/(\sqrt n-1)$.
These bounds count offsets. By \cref{lem:replication}, one offset $m$
produces $n-|m|$ intersecting pairs of original pieces. For example,
the even--odd decomposition in \cref{prop:imaginary} gives
$\A(i\sqrt n,n)=\{-1,0,1\}$, while $k(i\sqrt n,n)=n$ and
$J(i\sqrt n,n)=3n-2$. Bounded branching is therefore compatible with a
contact count that grows with the alphabet size.
\end{remark}

\section{First-digit loci and convex parameter geometry}\label{sec:piece-loci}

We now express the contact condition in parameter space. Fix $n\ge2$,
write $h=n-1$ and $N=2n-1$, and let $-h\le m\le h$ be an integer. Recall
from \eqref{eq:balanced} that $E(c,N)=K-K$ for $K=E(c,n)$.
The first-digit locus $\mathcal L_{n,m}$ consists of parameters for
which some overlap expansion \eqref{eq:offsets} begins with $m$.
Replacing the remaining digit choices by their convex hull gives the
convex locus $\widehat{\mathcal L}_{n,m}$:
\begin{align}
 \mathcal L_{n,m}
 &=\{c:|c|>1,\ 2c\in2m+c^{-1}E(c,N)\},\label{eq:piece-locus}\\
 \widehat{\mathcal L}_{n,m}
 &=\{c:|c|>1,\ 2c\in2m+c^{-1}\conv E(c,N)\}.
 \label{eq:convex-piece-locus}
\end{align}
Equivalently, $c\in\mathcal L_{n,m}$ if and only if
\begin{equation}\label{eq:piece-expansion}
 c=m+\sum_{j\ge1}d_jc^{-j},\qquad d_j\in[-h,h]\cap\Z.
\end{equation}
The convexification identity
\[
 \conv E(c,N)=\sum_{j\ge0}[-2h,2h]c^{-j}
\]
follows by taking convex hulls of finite digit sums and passing to the
Hausdorff limit; see also~\cite[Chapter~1, p.~42]{Bousch1992}.
Thus $\widehat{\mathcal L}_{n,m}$ is described by
\eqref{eq:piece-expansion} with real tail digits $d_j\in[-h,h]$,
while the first digit $m$ remains an integer. By \cref{lem:connectivity},
\begin{equation}\label{eq:piece-cover}
 \M_n=\bigcup_{m=-h}^{h}\mathcal L_{n,m},\qquad
 \A(c,n)=\{m:c\in\mathcal L_{n,m}\},\qquad
 \mathcal L_{n,m}\subseteq\widehat{\mathcal L}_{n,m}.
\end{equation}

\begin{figure}[tbp]
\centering
\includegraphics[width=\textwidth,alt={Eight quarter-plane panels for n = 7: the union of the first-digit loci, numbered by how many loci overlap, and the loci for m = 0 to 6 in blue; only the m = 6 locus stays inside the dashed circle of radius sqrt(7).}]{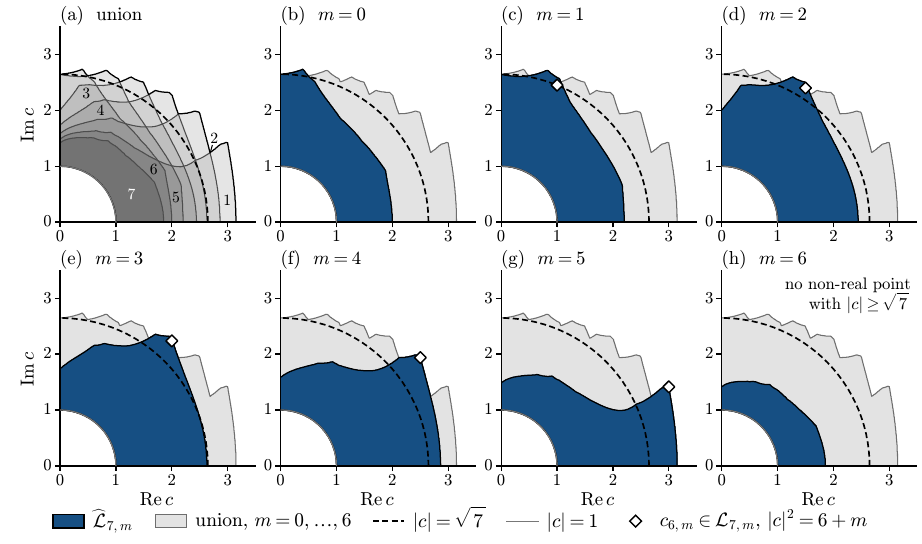}
\caption{Convex first-digit loci $\widehat{\mathcal L}_{7,m}$ for
$m=0,\ldots,6$ in the first quadrant, in blue against their union in
gray; panel (a) shows that union, numbered by how many loci contain
each point. The dashed arc is $|c|=\sqrt7$ and the thin gray arc is
$|c|=1$; white diamonds mark $c_{6,m}\in\mathcal L_{7,m}$ from
\eqref{eq:quadratic-label}, $1\le m\le5$. The loci
overlap near the unit disk, whereas the extreme locus $m=6$ has no
non-real part outside the dashed arc. These are numerical renderings of
the convex relaxation.}
\label{fig:piece-cells}
\end{figure}

These sets generally overlap. A disjoint decomposition is obtained by
recording the full set $\A(c,n)$ of admissible labels: for each nonempty
$S\subseteq[-h,h]\cap\Z$, take the parameters belonging to precisely
the loci labeled by $S$.
Complex conjugation preserves each label, while $c\mapsto-c$ exchanges
$m$ and $-m$. Hence there are $n$ representative labels up to symmetry
in the whole plane; negative labels may still occur in the first quadrant.
\Cref{fig:piece-cells} displays the nonnegative convex loci separately.

\subsection{Support functions and localization}

The convex relaxation has an explicit support description. Put
\[
 Z(c)=\sum_{j\ge1}[-1,1]c^{-j},\qquad
 \sigma_c(\phi)=\sum_{j\ge1}\lvert\operatorname{Re}(e^{-i\phi}c^{-j})\rvert.
\]
The sums converge for $|c|>1$, and $\sigma_c$ is the support function of
the centrally symmetric compact convex set $Z(c)$. Consequently
\begin{equation}\label{eq:support-membership}
 c\in\widehat{\mathcal L}_{n,m}
 \quad\Longleftrightarrow\quad
 \lvert\operatorname{Re}(e^{-i\phi}(c-m))\rvert\le h\sigma_c(\phi)
 \quad\text{for every }\phi\in[0,\pi].
\end{equation}
At a fixed parameter, the eligible real first digits form an interval:
they are the horizontal section of $c-hZ(c)$ on the real axis.
The convex-admissible integer labels therefore form a consecutive
interval, possibly empty. No such interval assertion is needed for the
exact labels.

For $n\ge3$, the radial enclosure for $\M_n\setminus\R$, and for the
larger set of expansions with real digits, is a theorem of Beaucoup,
Borwein, Boyd, and Pinner~\cite[Theorem~2]{BBBP1998}, stated there for
coefficient bounds $g>1$; it is applied to these loci in
~\cite[Proposition~4]{EspiguleJuherSaldana2024}. An elementary proof by
vertical support, valid for every $n\ge2$, is also given in
~\cite[Proposition~A.2]{EspiguleJuher2026Layers}. That elementary
estimate applies verbatim to the convex relaxation and also localizes
each first-digit locus in a vertical strip.

\begin{proposition}[Radial and strip bounds]\label{prop:piece-strips}
If $c=a+ib\in\widehat{\mathcal L}_{n,m}\setminus\R$ and $r=|c|$, then
\begin{equation}\label{eq:piece-strips}
 r<1+\sqrt h,\qquad |2a-m|\le\frac h{(r-1)^2}.
\end{equation}
In particular, the non-real part of $\widehat{\mathcal L}_{n,m}$ with
$r^2\ge n$ lies in the vertical strip
\[
 |2a-m|\le1+\frac2{\sqrt n-1}.
\]
\end{proposition}

\begin{proof}
Write $c=re^{i\theta}$ and use the real-tail expansion in
\eqref{eq:piece-expansion}. The estimate
$|\sin(j\theta)|\le j|\sin\theta|$ gives
\[
 |b|\le h\sum_{j\ge1}\lvert\operatorname{Im}(c^{-j})\rvert
 \le h|b|\sum_{j\ge1}j r^{-j-1}
 =\frac{h|b|}{(r-1)^2}.
\]
The second inequality is strict at a non-real parameter, because its
$j=2$ sine estimate is strict. Dividing by $|b|$ proves the radial
bound. Multiplying the expansion by $c$ gives
$c^2-mc=d_1+\sum_{j\ge1}d_{j+1}c^{-j}$. Its imaginary part gives
$|2a-m|\le h/(r-1)^2$. Finally,
$h/(\sqrt n-1)^2=1+2/(\sqrt n-1)$.
\end{proof}

For numerical rendering, truncate $Z(c)$ to
$Z_J(c)=\sum_{j=1}^{J}[-1,1]c^{-j}$. The geometric tail gives
\begin{equation}\label{eq:zonotope-tail}
 Z_J(c)\subseteq Z(c)\subseteq
 Z_J(c)+\overline B\left(0,\frac{|c|^{-J}}{|c|-1}\right).
\end{equation}
Finite zonotopes and their enlarged supporting half-planes thus give
inner and outer membership bounds for the convex relaxation. The
figures use this construction on sampled parameter grids; the tail
bound is explicit, but the floating-point renderings are not interval
certificates. The source code and grid data are among the ancillary
files described at the end of the paper.

\subsection{First digits outside the connectedness disk}

The disk estimate in \cref{lem:connectivity} singles out the range
$|c|\ge\sqrt n$ for the following classification. For each fixed
first digit $m$, both $\mathcal L_{n,m}$ and
$\widehat{\mathcal L}_{n,m}$ increase with $n$. The restriction below
concerns the portions outside this growing disk and first digits whose
magnitude increases with the alphabet size. Since $c\mapsto\bar c$
preserves every label and $c\mapsto-c$ exchanges $m$ and $-m$, positive
labels in the first quadrant are representative; the last assertion of
the theorem records the symmetric form.

\begin{theorem}[First digits outside the connectedness disk]\label{thm:exterior-labels}
Let $h=n-1$. The positive integers $m$ for which
$\mathcal L_{n,m}$ contains a parameter with
$\operatorname{Re}c>0$, $\operatorname{Im}c>0$, and $|c|^2\ge n$ are precisely
\begin{equation}\label{eq:surviving-labels}
 1\le m\le\min\{h,\lceil2\sqrt h\rceil\}.
\end{equation}
The same classification holds for $\widehat{\mathcal L}_{n,m}$.
Moreover, if $n\ge7$ and $c\in\widehat{\mathcal L}_{n,m}\setminus\R$
satisfies $|c|^2\ge\lceil2\sqrt h\rceil+1$, in particular if $|c|^2\ge n$,
then $|m|\le\lceil2\sqrt h\rceil$.
For every $m$ in \eqref{eq:surviving-labels}, the exact locus
$\mathcal L_{n,m}$ contains the parameter
\begin{equation}\label{eq:quadratic-label}
 c_{h,m}=\frac{m+1+i\sqrt{4h-(m-1)^2}}2,
 \qquad |c_{h,m}|^2=h+m.
\end{equation}
Consequently, for every integer $m$ with $|m|\le h$, each of the loci
$\mathcal L_{n,m}$ and $\widehat{\mathcal L}_{n,m}$ contains a non-real
parameter with $|c|^2\ge n$ if and only if
$|m|\le\min\{h,\lceil2\sqrt h\rceil\}$.
\end{theorem}

\begin{proof}
For $h\ge6$, put $q=\lceil2\sqrt h\rceil+1\le h$.
If $c\in\widehat{\mathcal L}_{n,m}$ with $|m|\ge q$, its expansion makes
$c^{-1}$ a zero of a series as in \cref{lem:jensen}, with real tail
coefficients. When $|c|^2\ge q$ and $c$ is non-real, this contradicts
that lemma. Hence $|m|\le q-1=\lceil2\sqrt h\rceil$ for every non-real
$c\in\widehat{\mathcal L}_{n,m}$ with $|c|^2\ge q$; since $q\le h<n$,
this includes the range $|c|^2\ge n$.
For $h\le5$, all positive allowed labels satisfy
$m\le\min\{h,\lceil2\sqrt h\rceil\}$.

For every label in \eqref{eq:surviving-labels}, one has
$m-1<2\sqrt h$, so \eqref{eq:quadratic-label} is non-real.
Its equation ${c^2-(m+1)c+h+m=0}$ gives
\begin{equation}\label{eq:constant-tail}
 c=m-\frac h{c-1}=m+\sum_{j\ge1}(-h)c^{-j}.
\end{equation}
These integer tail digits are admissible. Thus
$c_{h,m}\in\mathcal L_{n,m}$, and $|c_{h,m}|^2=h+m\ge n$.

For the final assertion, the exclusion above gives
$|m|\le\lceil2\sqrt h\rceil$ whenever $n\ge7$, while
$\min\{h,\lceil2\sqrt h\rceil\}=h$ for $n\le6$. Conversely, $c_{h,m}$
serves for $m>0$; for $m<0$, the map $c\mapsto-c$ carries $c_{h,|m|}$
into $\mathcal L_{n,m}$; and $i\sqrt n\in\mathcal L_{n,0}$, because
the proof of \cref{prop:imaginary} with $r^2=n$, so that $\eta=r$, gives
$P_0\cap(2i\sqrt n+P_0)=[a_0-1,a_0+1]+i\sqrt n\ne\varnothing$, whence
$0\in\A(i\sqrt n,n)$ by \cref{lem:replication}.
Since $\mathcal L_{n,m}\subseteq\widehat{\mathcal L}_{n,m}$, the same
parameters serve for the convex loci.
\end{proof}

\begin{proof}[Proof of \cref{thm:roots}]
Put $h=n-1\ge6$ and $q=\lceil2\sqrt h\rceil+1\le h$. The point
$\lambda=1/c$ is a non-real zero of
$F(z)=1+a_1z+\cdots+a_kz^k$, and $|\lambda|^2\le1/q$. If
$|a_1|\ge q$, then \cref{lem:jensen}, applied with $m=-a_1$ and
$b_j=-a_j$, shows that every zero of $F$ in $|z|\le q^{-1/2}$ is real,
a contradiction. Hence $|a_1|\le\lceil2\sqrt h\rceil$.

For the constructions, let $h\ge1$ and let $m$ be an integer with
$2\le m\le\min\{h,\lceil2\sqrt h\rceil\}$. By \eqref{eq:constant-tail},
$\lambda_0=1/c_{h,m}$ is a non-real zero of
$G(z)=1-mz+h\sum_{j\ge2}z^j$, which converges in the unit disk. The
partial sums $G_k$ of $G$ converge to $G$ locally uniformly there, so
Hurwitz's theorem gives zeros $\lambda_k\to\lambda_0$ of $G_k$. The
polynomials $z^kG_k(1/z)=z^k-mz^{k-1}+h(z^{k-2}+\cdots+1)$ are monic of
height $h$, with subleading coefficient $-m$, and have the roots
$1/\lambda_k\to c_{h,m}$. Since $c_{h,m}$ is non-real and
$|c_{h,m}|^2=h+m>n$, the same holds for $1/\lambda_k$ when $k$ is large.
Taking $m=\lceil2\sqrt h\rceil$ for $n\ge7$ proves sharpness, and taking
$m=h$ for $3\le n\le6$ gives non-real roots with $|c|^2>n$ and $|a_1|=n-1$.
\end{proof}

\begin{figure}[htbp]
\centering
\includegraphics[width=\textwidth,alt={Six quarter-plane panels for n = 3, 5, 6, 7, 12 and 21 with numbered first-digit bands; from n = 7 on, the blue bands of large first digits stay inside the dashed circle of radius sqrt(n).}]{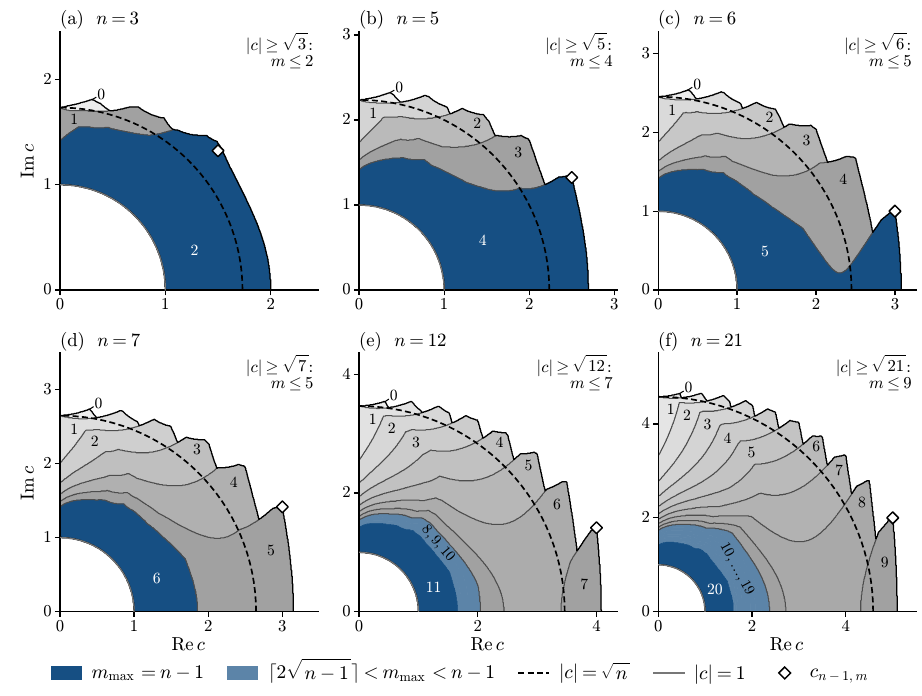}
\caption{First-quadrant convex first-digit geometry at six alphabet
sizes. The bands are level sets of the largest convex-admissible
nonnegative label $m_{\max}$, whose value is printed in each band; the labeled loci
themselves overlap. Dark blue marks $m_{\max}=n-1$ and lighter blue
$\lceil2\sqrt{n-1}\rceil<m_{\max}<n-1$; dashed arcs are $|c|=\sqrt n$
and thin gray arcs $|c|=1$. The upper right corner of each panel gives
the largest convex-admissible label at non-real grid points with
$|c|\ge\sqrt n$, namely
$\min\{n-1,\lceil2\sqrt{n-1}\rceil\}$, as in \cref{thm:exterior-labels};
white diamonds mark \eqref{eq:quadratic-label} with this $m$. Blue
extends outside $|c|=\sqrt n$ only for $n\le6$.}
\label{fig:piece-atlas}
\end{figure}

The extreme locus $\mathcal L_{n,n-1}$ therefore meets the non-real
range $|c|\ge\sqrt n$ exactly when $2\le n\le6$.
For $n\ge7$, the number of positive first digits whose convex loci
have no non-real parameter in that range is
\[
 n-1-\lceil2\sqrt{n-1}\rceil.
\]
This number is nondecreasing; it equals $1$, $4$, and $11$ at
$n=7$, $12$, and $21$, respectively. \Cref{fig:piece-atlas} shows
the transition and the subsequent loss of outer first-digit loci. For $3\le n\le6$,
the extreme-label parameter in \eqref{eq:quadratic-label} is exactly
the quadratic dendrite example in \eqref{eq:examples}.

\Needspace{6\baselineskip}
\begin{proposition}[Attainment of the contact bound]\label{prop:attainment}
For $h\ge4$, put $p=\lceil2\sqrt h\rceil$ and $n=h+1$.
At the parameter $c_h=c_{h,p}$ in \eqref{eq:quadratic-label},
\begin{equation}\label{eq:exact-minimum}
 k(c_h,n)=n-p.
\end{equation}
In particular, for every $n\ge7$,
\[
 \min_{c\in\M_n\setminus\R}k(c,n)
 =n-\lceil2\sqrt{n-1}\rceil.
\]
\end{proposition}

\begin{proof}
Since $p\le h$, \eqref{eq:constant-tail} shows that $p$ is an
admissible label. Also $2\operatorname{Re}c_h=p+1$ and
$|c_h|^2=h+p$. The strict estimate
\[
 \frac h{(\sqrt{h+p}-1)^2}<2
\]
is equivalent to $h/2+p-1-\sqrt{2h}>0$, which follows from
$p\ge2\sqrt h$ and $h\ge4$. By \cref{prop:piece-strips}, every
admissible integer label $m$ therefore satisfies $m>p-1$, hence $m\ge p$.
For a positive label $m$, the participating left indices are exactly
$\{m,\ldots,n-1\}$ by \eqref{eq:pair-offset}. These intervals are
nested, and $p$ itself occurs. Their union has $n-p$ indices.
Combining this equality with \cref{thm:multiplicity} proves the minimum.
\end{proof}

\begin{figure}[htbp]
\centering
\includegraphics[width=\textwidth,alt={The attractor K and its translate 2c + K as two rows of seven slanted pieces; P5 meets 2c + P0 at c - 1, P6 meets 2c + P1 at c + 1, and P6 does not meet 2c + P0.}]{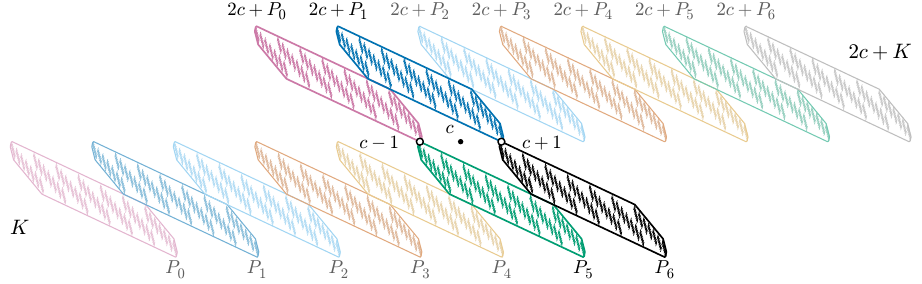}
\caption{The attractor $K=E(c,7)$ and its translate $2c+K$ at the
parameter $c=3+i\sqrt2$ of \cref{prop:attainment}, marked by the white
diamond in \cref{fig:piece-atlas}(d), with the convex hulls of their
first-level pieces; pieces with the same index share a hue.
The hulls of $P_i$ and $2c+P_j$ intersect exactly
when $i-j=5$, because $5$ is the only convex-admissible label at $c$
(\cref{prop:piece-strips,thm:exterior-labels}). In particular, the
rightmost piece $P_6$ of $K$ does not meet the leftmost piece $2c+P_0$
of the translate. By \eqref{eq:constant-tail},
$c-1\in P_5\cap(2c+P_0)$ and $c+1\in P_6\cap(2c+P_1)$ (white circles),
so two pieces of $K$ meet $2c+K$; these four pieces are drawn in full
color, the others in light tints. The symmetry $z\mapsto2c-z$ about $c$
(black dot) exchanges $K$ and $2c+K$, maps $P_i$ onto $2c+P_{6-i}$, and
swaps $c-1$ and $c+1$. The point samples are numerical illustrations.}
\label{fig:contact}
\end{figure}

\Cref{fig:contact} illustrates \cref{prop:attainment} for $h=6$, at
$c_6=3+i\sqrt2$, the white diamond in \cref{fig:piece-atlas}(d).

\suppressfloats[t]

\subsection{Confinement and rescaled geometry}

All first-digit loci retain a common part near the unit disk. For $c=ir$ with
$1<r^2\le n$, splitting the balanced digits into even and odd positions
as in \cref{prop:imaginary} gives
\[
 \A(ir,n)=\left\{m\in\Z:|m|\le
 \min\left\{h,\left\lfloor\frac h{r^2-1}\right\rfloor\right\}\right\}.
\]
Indeed, since $r^2\le n\le2n-1$, the digits in even and in odd positions
each fill an interval, so the tails $\sum_{j\ge1}d_j(ir)^{-j}$ fill the
rectangle $hZ(ir)$ displayed below, and $m\in\A(ir,n)$ exactly when
$ir-m$ lies in it.
In particular, every label is exactly admissible when $1<r^2\le2$.
For $1<r^2<2$, the point $ir-m$ lies strictly inside the tail rectangle
\[
 hZ(ir)=\left[-\frac h{r^2-1},\frac h{r^2-1}\right]
 +i\left[-\frac{hr}{r^2-1},\frac{hr}{r^2-1}\right]
\]
for every label $m$. Hausdorff continuity of $c\mapsto Z(c)$, which
follows from \eqref{eq:zonotope-tail}, then gives a common open
intersection of all convex loci, extending into the first quadrant.

First digits comparable to the alphabet size satisfy a stronger
radial bound.

\begin{proposition}[Confinement of large first digits]\label{prop:edge-confinement}
Let $h=n-1\ge6$, and let $m$ be an integer with
$1+2\sqrt h\le m\le h$. Set
\begin{equation}\label{eq:edge-radius}
 R_h(m)=\frac{2\sqrt{1+2h^2}}
 {(2/3)^{3/2}\bigl(m+1+\sqrt{(m-1)^2-4h}\bigr)}.
\end{equation}
Every non-real $c\in\widehat{\mathcal L}_{n,m}$ satisfies
$|c|^2\le R_h(m)$. For any sequence of integers $m_h\in[1,h]$
such that $m_h/h\to\tau\in(0,1]$, the bound applies for all sufficiently
large $h$, and
\begin{equation}\label{eq:edge-collapse}
 R_h(m_h)\longrightarrow\frac{3\sqrt3}{2\tau}.
\end{equation}
\end{proposition}

\begin{proof}
The majorant in \eqref{eq:majorant-root} provides a real zero
$\alpha\le t_0$. If $|c|^2>3/2$, the zeros $\alpha$, $c^{-1}$, and
$\bar c^{-1}$ lie in $|z|<\rho=\sqrt{2/3}$, and Jensen and Parseval give
$\rho^3|c|^2/\alpha\le(1+\tfrac23m^2+\tfrac43h^2)^{1/2}\le\sqrt{1+2h^2}$,
yielding the bound.
The remaining case is automatic because $t_0>1/m\ge1/h$ makes
$R_h(m)>3\sqrt3/2>3/2$. If $m_h/h\to\tau>0$, then
$m_h\ge1+2\sqrt h$ for all sufficiently large $h$.
Dividing the numerator and denominator in \eqref{eq:edge-radius} by
$h$ gives \eqref{eq:edge-collapse}.
\end{proof}

Thus, when $m_h/h\to\tau>0$, the non-real portions of the convex loci
remain in a bounded disk while the connectedness disk has radius
$\sqrt n$. Their distances from the origin, after division by $\sqrt n$,
tend uniformly to zero. The constant in \eqref{eq:edge-collapse} is not
optimal; a choice of $\rho$ depending on $\tau$ improves it. For
$\tau<0.37$, \cref{prop:piece-strips} alone gives the better asymptotic
bound $|c|\le1+\tau^{-1/2}+o(1)$, because
$m_h\le2|c|+h/(|c|-1)^2$ and $|c|<1+\sqrt h$.

\begin{figure}[htbp]
\centering
\includegraphics[width=\textwidth,alt={Three quarter-plane panels for n = 7, 21 and 100 in the coordinate c divided by sqrt(n - 1): toothed first-digit loci in an annulus around the unit circle that narrows as n grows, darker toward the real axis.}]{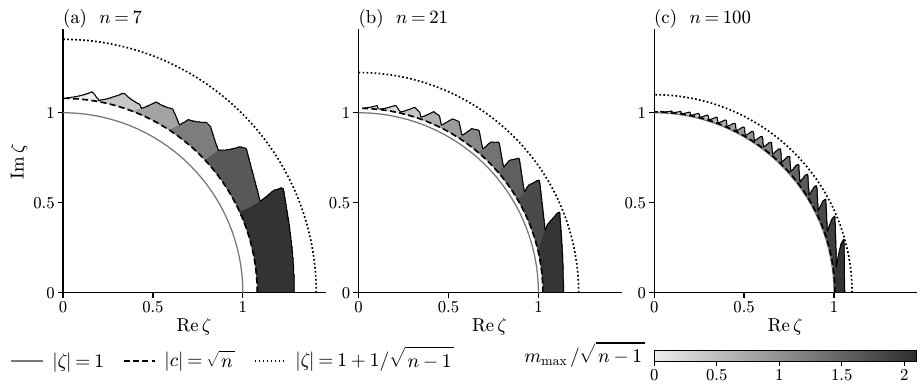}
\caption{Convex first-digit loci outside the connectedness disk for
$n=7$, $21$, and~$100$, in the coordinate $\zeta=c/\sqrt{n-1}$, shaded by
$m_{\max}/\sqrt{n-1}$ (as in \cref{fig:piece-atlas}) on a common scale.
The thin gray arc is the unit circle; the dashed arc $|c|=\sqrt n$ and
the dotted arc $|\zeta|=1+1/\sqrt{n-1}$ bound the band of
\eqref{eq:normalized-bounds}. The localization
\eqref{eq:normalized-geometry} explains the narrowing band and the
relation between the first digit and the horizontal position.}
\label{fig:normalized-pieces}
\end{figure}

For the portions outside the connectedness disk, shown in
\cref{fig:normalized-pieces}, set $\zeta=c/\sqrt h$.
If $c\in\widehat{\mathcal L}_{n,m}\setminus\R$ and $|c|^2\ge h+1$,
then \cref{prop:piece-strips} gives the explicit inequalities
\begin{equation}\label{eq:normalized-bounds}
 \sqrt{1+\frac1h}\le|\zeta|<1+\frac1{\sqrt h},\qquad
 \left|\frac m{\sqrt h}-2\operatorname{Re}\zeta\right|
 \le\frac{1+2/(\sqrt{h+1}-1)}{\sqrt h}.
\end{equation}
In particular, uniformly over these parameters and first digits,
\begin{equation}\label{eq:normalized-geometry}
 |\zeta|=1+O(h^{-1/2}),\qquad
 \frac m{\sqrt h}=2\operatorname{Re}\zeta+O(h^{-1/2}).
\end{equation}
The rescaled non-real outer portions therefore concentrate near the unit circle,
with their first digits ordered approximately by horizontal position.
The estimates describe their location
without determining the topology of the individual loci.

\FloatBarrier

\section{Examples at every alphabet size from two to six}
\label{sec:examples}

We complete the proof of the sharp bound by exhibiting a dendrite at each
of the five admissible alphabet sizes. Up to affine normalization and
complex conjugation, the binary example is the one in
~\cite[Figure~3]{ErogluRohdeSolomyak2010}. The calculations below give
a self-contained verification of all five contact structures, using
integer arithmetic and one rational interval for the binary case.

\begin{proposition}[Dendrites for $2\le n\le6$]\label{prop:all-five-examples}
For each \(2\le n\le6\), there is a non-real parameter \(c\) for which
\(E(c,n)\) is a dendrite. One may take the upper-half-plane root of
\[
 \begin{cases}
 c^3-2c^2+2=0,&n=2,\\
 c^2-nc+2(n-1)=0,&3\le n\le6.
 \end{cases}
\]
At each of these parameters,
\[
 E(c,n)\cap\bigl(2c+E(c,n)\bigr)=\{c\},
\]
and non-neighboring first-level pieces are disjoint. Thus the first-level
intersection graph is a path, with neighboring pieces meeting at the
points \(t+1\), \(t\in\{-n+1,-n+3,\ldots,n-3\}\).
\end{proposition}

\begin{proof}
\emph{Difference expansions and inverse states.}
Write \(h=n-1\), \(D_n=\{-h,\ldots,h\}\), and \(K=E(c,n)\).
An intersection between pieces separated by \(2j\) gives an expansion
\[
 jc=\sum_{\ell\ge0}e_\ell c^{-\ell},\qquad e_\ell\in D_n.
\]
Conversely, every such expansion can be lifted coordinatewise to a pair
of addresses in \(A_n\), and hence gives an intersection. Its inverse
states satisfy
\[
 z_0=jc,\qquad z_{\ell+1}=c(z_\ell-e_\ell),\qquad
 |z_\ell|\le R:=\frac{h}{1-|c|^{-1}}.
\]
We enumerate the states reachable from \(jc\) while satisfying a
strict upper bound on the squared norm; these states play the role of the
neighbor maps of Bandt and Graf~\cite{BandtGraf1992}. In the five cases below, the
reachable graphs are finite. Repeatedly deleting vertices with no
outgoing edge identifies all infinite admissible paths.

\smallskip
\emph{The quadratic examples.}
Suppose \(3\le n\le6\). Represent a state as
\(z=a+bc\), where \(a,b\in\mathbb Z\). Since
\(c+\bar c=n\) and \(c\bar c=2h\), its transition and squared norm are
\begin{equation}\label{eq:quadratic-certificate}
 (a,b)\xrightarrow{e}(-2hb,a+nb-e),\qquad
 H_n(a,b)=a^2+nab+2hb^2.
\end{equation}
For \(n=3,4,5,6\), respectively, choose \(B_n=17,26,39,54\).
Then \(R^2<B_n\). Indeed, if
\(w_n=B_n(2h+1)-2h^3\), the required inequality follows from
\(w_n>0\) and \(w_n^2>8hB_n^2\), both direct integer inequalities.
Thus every bounded inverse path has \(H_n(a,b)<B_n\).

The following table lists every transition reachable from
\(c=(0,1)\) and from \(2c=(0,2)\) that satisfies this norm bound. Entries \(e:(a,b)\) specify the digit and
successor. A dash means that no successor satisfies the bound.
\begin{center}
\small
\begin{tabular}{ccl}
\toprule
\(n\)&State \((a,b)\)&Successors with \(H_n<B_n\)\\
\midrule
3&\((0,1)\)&\(0:(-4,3),\ 1:(-4,2),\ 2:(-4,1)\)\\
 &\((-4,1)\)&\(-2:(-4,1),\ -1:(-4,0)\)\\
 &\((-4,0)\)&\(-2:(0,-2)\)\\
 &\((-4,2),(-4,3),(0,-2),(0,2)\)&---\\
\midrule
4&\((0,1)\)&\(1:(-6,3),\ 2:(-6,2),\ 3:(-6,1)\)\\
 &\((-6,1)\)&\(-3:(-6,1)\)\\
 &\((-6,2),(-6,3),(0,2)\)&---\\
\midrule
5&\((0,1)\)&\(1:(-8,4),\ 2:(-8,3),\ 3:(-8,2),\ 4:(-8,1)\)\\
 &\((-8,1)\)&\(-4:(-8,1)\)\\
 &\((-8,2),(-8,3),(-8,4),(0,2)\)&---\\
\midrule
6&\((0,1)\)&\(1:(-10,5),\ 2:(-10,4),\ 3:(-10,3),\)\\
 &&\(4:(-10,2),\ 5:(-10,1)\)\\
 &\((-10,1)\)&\(-5:(-10,1)\)\\
 &\((-10,2)\)&\(-5:(-20,7),\ -4:(-20,6),\ -3:(-20,5)\)\\
 &\((0,2)\)&\(5:(-20,7)\)\\
 &\((-10,3),(-10,4),(-10,5)\)&---\\
 &\((-20,5),(-20,6),(-20,7)\)&---\\
\bottomrule
\end{tabular}
\end{center}
After repeatedly deleting vertices with no outgoing edge, the unique
infinite path is
\[
 c\xrightarrow{h}c-2h\xrightarrow{-h}c-2h\xrightarrow{-h}\cdots.
\]
The state \(2c\) has no successor satisfying the bound for \(n=3,4,5\);
for \(n=6\), its only such successor is \((-20,7)\), which has no
successor satisfying the bound.
For every \(j\ge3\), \(H_n(0,j)=2hj^2\ge B_n\).
It follows that \(c\) has the unique expansion
\[
 c=h-hc^{-1}-hc^{-2}-\cdots,
\]
whereas \(jc\) has none for \(2\le j\le n-1\).

\smallskip
\emph{The binary example.}
For \(n=2\), both critical values of
\(P(x)=x^3-2x^2+2\) are positive, so \(P\) has a unique real root
\(\beta\). Direct rational evaluation gives
\(-84/100<\beta<-83/100\).
Its other roots are \(c\) and \(\bar c\), and
\(c\bar c=-2/\beta>9/4\). Hence \(|c|>3/2\) and \(R<3\).
Represent a state as \(a+bc+dc^2\). Its transition is
\[
 (a,b,d)\xrightarrow{e}(-2d,a-e,b+2d).
\]
Its squared modulus is \(v_0+v_1\beta+v_2\beta^2\), where
\begin{align*}
 v_0&=a^2+2ab+4ad+2bd+4d^2,\\
 v_1&=-ab-2b^2-4bd-2d^2,\\
 v_2&=-ad+b^2+2bd.
\end{align*}
Substitution of the rational interval for \(\beta\) gives the following
complete list of successors with squared modulus less than \(9\).
Every excluded successor has squared modulus strictly greater than \(9\).
\begin{center}
\small
\begin{tabular}{cl}
\toprule
State \((a,b,d)\)&Successors with squared modulus less than \(9\)\\
\midrule
\((0,1,0)\)&\(0:(0,0,1),\ 1:(0,-1,1)\)\\
\((0,0,1)\)&\(1:(-2,-1,2)\)\\
\((0,-1,1)\)&\(-1:(-2,1,1),\ 0:(-2,0,1),\ 1:(-2,-1,1)\)\\
\((-2,-1,1)\)&\(-1:(-2,-1,1)\)\\
\((-2,0,1)\)&\(-1:(-2,-1,2),\ 0:(-2,-2,2)\)\\
\((-2,1,1),(-2,-1,2),(-2,-2,2)\)&---\\
\bottomrule
\end{tabular}
\end{center}
Thus the unique infinite path from \(c\) has digits
\(1,1,-1,-1,\ldots\), giving
\[
 c=1+c^{-1}-c^{-2}-c^{-3}-\cdots.
\]
The displayed identity also follows directly from \(P(c)=0\).

\smallskip
\emph{Contact structure and dendricity.}
In all five cases the unique difference expansion uses only extreme
digits \(\pm h\). Each such digit has a unique lift to a pair from
\(A_n\), namely \((h,-h)\) or \((-h,h)\). Therefore the neighboring
contact is exactly \(\{c\}\). The absence of expansions of \(jc\)
for \(j\ge2\) proves that non-neighboring pieces are disjoint.

Neighboring contact gives connectedness of \(K\) by
\cref{lem:connectivity}, and Hata's theorem gives local connectedness
~\cite{Hata1985}. The first-level intersection graph is a path. Every
first-level junction is a cut point separating the pieces on its two sides.
A simple closed curve cannot cross such a junction, since a circle
minus one point remains connected. Consequently every simple closed
curve in \(K\) would lie in a single first-level piece. Repeating inside
that piece traps the curve in nested cylinders of diameter tending to
zero, a contradiction. Hence \(K\) is a dendrite; compare the
criterion of Bandt and Keller~\cite{BandtKeller1991}.
\end{proof}

\section{Further questions}

For $2\le n\le6$, it remains to classify the dendrite parameters and
their possible contact sets. With $h=n-1$, \cref{prop:one-piece} locates
every non-real dendrite parameter in
\[
 (\mathcal L_{n,h}\cup\mathcal L_{n,-h})\setminus
 \bigcup_{|m|<h}\mathcal L_{n,m},
\]
with the first-digit loci of \cref{sec:piece-loci}. The condition
$k(c,n)=1$ counts participating pieces; it does not assert that the
contact is a singleton. For $n=2$, a singleton contact makes a connected
attractor a dendrite, by the argument at the end of \cref{sec:examples}
(compare~\cite{BandtKeller1991}). The converse fails, so the equivalence
stated in~\cite[Section~1]{ErogluRohdeSolomyak2010} holds only in the
forward direction. For real $c$ with $1<|c|<2$, the attractor $E(c,2)$
is an arc whose two pieces overlap in an interval. For a non-real
example, take $c=1/s\approx0.9334-1.1325i$, where
$s\approx0.4334+0.5258i$ is the root of $z^4+z^3+z^2-z+1$ in
~\cite[Example~3.8]{CalegariWalker2026}. The equation for $s$ is
equivalent to $f_{-1}f_1f_{-1}f_{-1}f_{-1}=f_1f_{-1}f_1f_1f_1$, so the
two first-level pieces share a fifth-level piece. Calegari and Walker
note that this piece is their entire intersection and that the attractor
is a dendrite. It would be useful to determine which
non-singleton contacts occur at the other admissible alphabet sizes.

For $n\ge7$, the exact minimum is known, but its full set of attaining
parameters remains to be described. The first-digit loci also raise
questions about the topology of their incidence strata and the difference
between exact membership and convex relaxation. The quadratic parameters
and support estimates in \cref{sec:piece-loci} give explicit starting points.

The proofs use both the common linear part and the equal spacing of the
fixed points; equal spacing drives the replication of contacts in
\cref{lem:replication}. It would be interesting to know whether collinear
systems with a common non-real multiplier but unevenly spaced fixed
points can produce dendrites with seven or more maps.

\section*{Acknowledgments}
The author thanks David Juher and Joan Salda\~na for discussions on
collinear self-similar sets.
The research program was supported by project PID2023-146424NB-I00,
funded by MICIU/AEI/10.13039/501100011033 and ERDF/EU, and by the
Universitat de Girona and Banco Santander Grant Programme for
Researchers in Training (IFUdG 2022--2024).

\section*{Computational materials}
The exact-arithmetic verifier for the tables in \cref{sec:examples}
(which uses only the Python standard library), the figure programs, and
the sampled parameter data are available as ancillary files with the
arXiv version of this article.
The Collinear Fractals explorer~\cite{Espigule2026Software} provides
interactive views of the attractors and connectedness loci:\\
\url{https://complextrees.com/collinear-fractals-gpu/}.

\section*{Use of artificial intelligence}
ChatGPT (OpenAI) assisted in developing and checking the proofs in
\cref{sec:topology,sec:jensen,sec:piece-loci}, preparing the verification
and visualization code, and drafting and revising the manuscript.
Claude (Anthropic) assisted in checking the proofs, computations, and
references, in drafting and revising parts of the text, including
\cref{thm:roots} and \cref{rem:elementary}, and in preparing
\cref{fig:contact}. The author is responsible for the mathematical
content.

\end{document}